\documentclass[11pt]{amsart}
\usepackage{amsfonts,color,graphicx}
\usepackage{eucal}
\usepackage{amsmath,amssymb,amsfonts,amsthm,enumerate}
\usepackage{mathrsfs}

\newtheorem{theorem} {{ Theorem}}[section]
\newtheorem{proposition} {{ Proposition}}[section]

\newtheorem{lemma} {{ Lemma}}[section]
\newtheorem{corollary} {{ Corollary}}[section]

\newtheorem{remark} {{ Remark}}[section]

\newtheorem{definition} {{ Definition}}[section]

\def\R{{I\!\!R}}

\newcommand{\be} {\begin{equation}}
\newcommand{\ee} {\end{equation}}
\newcommand{\bea} {\begin{eqnarray}}
\newcommand{\eea} {\end{eqnarray}}
\newcommand{\Bea} {\begin{eqnarray*}}
\newcommand{\Eea} {\end{eqnarray*}}
\newcommand{\pa} {\partial}

\newcommand{\Om} {\Omega}

\newcommand{\la} {\lambda}

\newcommand{\vp} {\varphi}
\newcommand{\var} {\varepsilon}
\newcommand{\ml} {{\mathscr L}}
\newcommand{\mh} {{\mathscr H}}
\newcommand{\supp}{\text{supp }}

\title{Two- and Multi-phase Quadrature Surfaces}

\author{Avetik Arakelyan}
\address{Institute of Mathematics, National Academy of Sciences of Armenia, 0019 Yerevan, Armenia}
\email{arakelyanavetik@gmail.com}

\author{Jyotshana V. Prajapat}
\address{Department of Mathematics, University of Mumbai, Vidyanagari, Santacruz (east), 400 097 Mumbai, India}
\email{jvprajapat@gmail.com}

\author{ Henrik Shahgholian}
\address{Department of Mathematics,  Royal Institute  of Technology,
100 44 Stockholm, Sweden}
\email{henriksh@math.kth.se}

\thanks{H.~ Shahgholian is partially supported by the Swedish Research
  Council}


\begin{document}

\begin{abstract}

In this paper we shall initiate the study of the two- and multi-phase quadrature surfaces (QS), which amounts to a two/multi-phase free boundary problems of Bernoulli type.
The problem is studied mostly from a potential theoretic point of view that (for two-phase case) relates to integral representation
$$
\int_{\partial \Omega^+} g h (x) \  d\sigma_x - \int_{\partial \Omega^-} g h (x) \  d\sigma_x= \int h d\mu \ ,
$$
where $d\sigma_x$ is the surface measure, $\mu= \mu^+ - \mu^-$ is given measure  with support in (a priori unknown domain) $\Omega$, 
$g$ is a given smooth positive function,
 and the integral holds for all  functions $h$, which are  harmonic on $\overline \Omega$.

Our approach is based on   minimization of the corresponding two- and multi-phase  functional and the use of its one-phase version as a barrier. We prove  several results concerning existence, qualitative behavior, and regularity theory for solutions. A central result in our study states  that  three or more junction points do not appear. 
\end{abstract}

\maketitle
\tableofcontents

\section {Introduction}
The current paper concerns the so-called quadrature identities for surface integrals, for the harmonic class of functions, and for given measures. Our primary purpose is to generalize the concept of quadrature surface (henceforth QS) to the two- and multi-phase counterpart.

The free boundary problem studied here has some "new" components that might be interesting to 
free boundary and  potential theory community. From potential theory point of view, we consider here a completely new problem  dealing with  the two-phase version of the problem of gravi-equivalent bodies\footnote{Here one body is the given measure, and the second body is a thin shell. The latter is the boundary of a domain, containing the other body.}, in particular, the existence of surfaces that "surround" the body is  essential  hearth of matter. On the other hand, the free boundary communities, specially those working with regularity theory,  would find an interesting extension of the concept of  two-phase Bernoulli problem, with the zero set having non-void   interior. This obviously makes the problem a three phase problem with the third phase being free of fluid.

\subsection{One-phase QS}\label{1phase}
Let $\Om \subset \R^N$ ($N\geq 2$)  be a bounded domain with reasonably smooth boundary, and $\mu$ be a measure with 
 support contained in $\Om$. Then we say that $\partial \Om$ is a quadrature surface with respect to $\mu$ if 
 the overdetermined Cauchy problem
\be\label{od}\left\{
\begin{array}{lll}
 \Delta u = - \mu~~~~~ &{\rm in  }&~ \Om \\
u =0, ~\frac{\pa u}{\pa \nu} =-1~~~~~~ &{\rm on }&~ \pa \Om ,
\end{array} \right.
\ee
has a solution. Here $\nu$ is the outward normal to the boundary $\pa \Om$.

For a better understanding, we recall the definition of one phase quadrature domains from \cite{gs}:  Given density  functions $  0 \leq g $, $ h \in L^\infty(\R^N)$ and  a Radon measure $\mu$,  we say that $\Om$ is a {\em quadrature domain} for $\mu$, for the given densities $g$ and $h$ if $\Om$ is a bounded open set in $\R^N$ such that 
\bea
&& {\rm supp}~ \mu \subset \Om,\\
&& U^\nu = U^\mu \quad {\rm in } \quad \R^N \setminus \Om ,\\
&& {\rm where} \quad \nu =    h \ml^N\mid_{\Om} + g  \mh^{N-1}\mid_{\pa \Om} \label{nu},
\eea
denoting the $n$  dimensional Lebesgue measure as $\ml^N$ and the $N-1$ dimensional Hausdorff measure as ${\mathcal H}^{N-1}$.
Here 
\be
U^\mu (x) = \int G(x-y) d\mu(y), \, x \in \R^N \ee denotes the Newtonian potential  corresponding to the measure $\mu$ with
 \be
G(x)=\left\{
\begin{array}{ll}
     \frac{1}{(N-2) \omega_N |x|^{N-2}} & \quad {\rm for } \quad N \geq 3, \\
 \frac{1}{2 \pi} \log |x| & \quad {\rm for } \quad N = 2 \end{array}\right. , \ee 
and hence
\be
- \Delta U^\mu = \mu.
\ee 
Let ${\mathcal W}$  be the set of all harmonic functions which can be expressed as linear combinations of $\{ \eta_y(x) := G(x-y)\}_{ y \in \R^N}$. Then, it can be verified  that $\Om$ is a quadrature domain if and only if the {\em quadrature identity\/}
\be\label{qi}
\int\limits_\Om \eta \, d\mu = \int\limits_\Om \eta  h \, dx + \int\limits_{\pa\Om} \eta  g \, d \mh^{N-1}\quad {\rm for ~all ~} \eta \in {\mathcal W}.
\ee

Quadrature domains can be obtained as supports of  local minimizers for the one phase functional
\be\label{j1}
J^1_{f,g} =   \int\limits_{\R^N} \left( |\nabla u|^2 - 2f  u  + g^2 \chi_{\{u > 0\}} \right)  \, dx,
\ee
where $f$, $g \in L^\infty(\R^N)$ are suitably chosen and  satisfy suitable conditions to allow a minimum for the functional.
 It was shown in \cite{gs} that a local minimum of the functional $J^1_{f,g}$ satisfies
\be\label{i1}
\begin{cases}
	 \Delta u = -f & {\rm in  }\;\; \Om = \{ u > 0 \} ,\\
      u =0,\; |\nabla  u| = g & {\rm on } \;\;\pa \Om .
\end{cases}
\ee
For general measures, e.g. Dirac masses,  the functional may not have lower bound, and hence the minimization may  not work. However, there is an easy way out of this problem, by smoothing out the measure and solving the approximate problem, and then considering  the limit problem. Indeed, for a given measure $\mu$, one uses radial mollifiers,  $\tilde \mu$ for approximating $\mu$. For $f =  \tilde \mu - h $ let  $\tilde u\geq 0 $ denote local minimum for $J_{\{f, g\}}$ so that it satisfies the equation (\ref{i1}).   See \cite{gs} for  details.

 Equation  (\ref{i1})  can be rewritten in the  sense of distributions as
\be
\Delta u + f \ml^N\mid_{\Om} = g {\mh}^{N-1}\mid_{\pa \Om}, \quad \Om = \{ u > 0\}. \ee 
In terms of the measure  $\nu $ defined in  (\ref{nu}), the above  identity can be written as $\mu +\Delta u = \nu$, so that $u$ is the difference of Newtonian potentials for the measures $\mu$ and $\nu$. The set  $\Om$ is a quadrature domain for $\mu$ if and only if ${\rm supp} (\mu) \subset \Om$. Now if we let $h =0$ then $\nu$ is the surface measure, and solution to this problem corresponds to (one phase) {\em quadrature surfaces}.

\subsection{Two-phase model}
 The two-phase counterpart of  the  functional (\ref{j1})   is
\be\label{j2}
J_{\{f_1, f_2,g\}}(u) :=  \int\limits_{\R^N} \left\{ |\nabla u|^2 - 2f_1 u^+ + 2f_2 u^- + g^2 \chi_{\{u \neq 0\}} \right\}  \, dx
\ee 
for  given functions $f_1$, $f_2$, $g$, where $u^+(x) := \max \{ u(x), 0 \}$ and $u^-(x) := \max \{ - u(x), 0 \}$.
The two phase functional with $g=0$, i.e., $J_{\{f,g=0\}}$ was studied in  the paper \cite{bjh}.  In this paper we are interested in showing existence of  ``two-phase'' quadrature surfaces corresponding to a measure $\mu =  \mu^+ - \mu^-$.  Thus,  assume that   ${\rm supp} (\mu^+) \neq \emptyset$,  ${\rm supp}(\mu^-) \neq \emptyset$  and that  ${\rm supp} ( g^2 ) $  has positive measure. We look for minimizer of the functional (\ref{j2}) where the  functions $f_1$ correspond to mollification of $\mu^+$ and $f_2$  is mollification of $\mu^-$. Here we expect that the (local) minimizer of (\ref{j2}) will satisfy 

\begin{align}\label{twophase}
	\begin{cases}
		 \Delta u = -f_1 \chi_{\{u>0\}} + f_2 \chi_{\{u<0\}}\quad &{\rm in} ~\Om, \\
		 u = 0, |\nabla  u| =  g \quad &{\rm on } ~\pa \Om,	
	\end{cases}	
\end{align}
where 
\be \Omega =  {\rm int}( \overline{\Om^+\cup\Om^-}),
~~~\Om^\pm :=\{ x \in \R^N: \pm u(x) > 0\}.\ee

Our approach in proving existence of minimizers to the two-phase functional, shall follow that of \cite{bjh}. By  relating the two phase functional to the one phase functionals, one  can efficiently generate solutions to the two phase problems by using suitable conditions ensuring existence of one phase solutions.

As mentioned earlier,   our problem produces three different phases, rather than two. More exactly, and contrary to the classical Bernoulli-type free boundaries, the interior of the set $\{u=0\}$ is  non-void in our case. In particular, one has  a triple junction free boundary points, where all three phases meet.
This type of  Bernoulli-free boundary is subject of current  study  by the third author, and his collaborators, see \cite{ASW}.

\begin{remark}\label{rem-1}
The general case where one replaces $g^2\chi_{\{u\neq 0\}}$ with 
$g_1^2\chi_{\{u > 0\}}+ g_2^2\chi_{\{u < 0\}}$ is not treated in this paper, but 
can be handled in much the same way as our situation. The functional, in this general case leads to the Bernoulli condition
$$
|\nabla u^+|^2 - |\nabla u^-|^2= g_1^2 -g_2^2
$$
on the two-phase boundary, and the standard one-phase boundary condition holds on one-phase boundary points.
\end{remark}

\section{Notation}
Here, for the reader convenience, we present some notations, which will be used  during  the exposition of the paper:
\begin{align*}
H^1(D)&\quad\{u\in L^2(D) : \nabla u\in L^2(D)\}\\
HL^1(D)&\quad \mbox{set of all integrable harmonic functions over $D$}\\
N &\quad \mbox{space dimension}\\
c, c_i, C_N,\alpha_i &\quad \mbox{generic constants}\\
\mathscr{H}^N&\quad\mbox{N-dimensional Hausdorff measure} \\
\mathscr{L}^N&\quad \mbox{N-dimensional Lebesgue measure}\\
\mu,\nu &\quad \mbox{Radon measures}\\
supp(\mu) &\quad \mbox{support of $\mu$}\\
\chi_D &\quad \mbox{the characteristic function of the set $D$}\\
\overline{D} &\quad \mbox{the closure of the set $D$}\\
{\rm int}(D) &\quad \mbox{interior of $D$}\\
\partial D &\quad \mbox{the boundary of $D$}\\
B(x,\tau),B_\tau(x) &\quad \{y\in\R^N : |y-x|<\tau \}\\
|B_\tau| &\quad \mbox{volume of a ball $B_\tau(x)$ }\\
d\sigma_x  &\quad \mbox{surface measure}\\
\cfrac{\partial\psi}{\partial \nu} &\quad \mbox{outward normal derivative of a function $\psi$}\\
\phi*\psi &\quad \mbox{convolution of  $\phi$ and $\psi$}\\
\delta_x  &\quad \mbox{Dirac measure at $x\in \R^N$}\\
\oint\limits_{\partial D} u\;d\mathscr{H}^{N-1}  &\quad \mbox{ the average integral of $u$ over $\partial D$}\\
\mu\mid_{D}&\quad \mbox{the restriction of $\mu$ to the set $D$}\\
\end{align*}


\section{Existence of minimizers}
\setcounter{equation}{0}
In this section, we give some conditions for existence of minimizers of the variational functional (\ref{j2}). 
We begin by proving the following comparison lemma for $J_{\{f_1, f_2, g\}}$, similar to Lemma 1.1 in  \cite{fp,gs}. 

\begin{lemma}\label{l1}
Assume that $f_1 \leq \tilde f_1$, $f_2 \leq \tilde f_2$ and $g  \geq \tilde g$. Let us denote $J
:= J_{\{f_1, f_2, g\}}$ and $\tilde J
:= J_{\{\tilde f_1, \tilde f_2, \tilde g\}}$. For  every $u_1, u_2 \in H^1(\R^N)$, we have  $v=\min \{u_1, u_2\}\in H^1(\R^N), w=\max \{u_1, u_2 \}\in H^1(\R^N),$ and 
\be
J (v) + \tilde J(w) \leq J(u_1) + \tilde J(u_2). 
\ee

\end{lemma}
\proof
We use the result that for a nondecreasing function $\Phi : \R \to \R$ and functions $h_1$, $h_2$ such that $h_1 \leq h_2$ we have
\be\label{33}
\int \left( h_1 \Phi(z_1) +  h_2 \Phi(z_2) \right) dx \leq  \int \left(h_1 \Phi(\min\{z_1, z_2\}) +  h_2 \Phi(\max\{z_1, z_2\})\right) dx, 
\ee
for any integrable functions $z_1$ and $z_2$.
Letting $\Phi(t) = \max\{t,0\}$ and $z_i = u_i$, $i = 1,2$, $h_1 = f_1$, $h_2 = \tilde f_1$  we get
\be
 \int \left( f_1 u_1^+  + \tilde f_1 u_2^+ \right) \, dx \leq  
 \int  \left( f_1 v^+ +  \tilde f_1 w^+ \right) \, dx.
 \ee
Taking again $z_i =  u_i,$ $i = 1,2$, $h_1 = f_2$, $h_2 = \tilde f_2$  and $\Phi(t) = \min\{t,0\}$  we get
\be
 \int \left( f_2 (-u_1^-)  + \tilde f_2(- u_2^-) \right) dx \leq \int  \left( f_2(- v^-) +  \tilde f_2 (- w^-) \right) dx.
  \ee
Finally, if we let $ h_1 = -g^2$, $h_2 = - \tilde g^2$ and $\Phi$ be the Heaviside function, then from (\ref{33}) we have the inequality
\be
\int \left( (-g^2) \chi_{\{u_1\neq 0\}} + (-\tilde g^2) \chi_{\{u_2\neq 0\}}  \right)  dx \leq  \int \left( (-g^2) \chi_{\{v\neq 0\}} + (-\tilde g^2) \chi_{\{w\neq 0\}}  \right) dx.
\ee
Since $\int \left( |\nabla u_1|^2 +  |\nabla u_2|^2 \right) \, dx = \int \left( |\nabla w|^2 +  |\nabla v|^2 \right) \, dx$ then  we conclude our inequality.
\qed

Observe that $J\mid_{\{u \in H^1(\R^N): u \geq 0\}} = J^1_{\{f_1,g\}}$ is a one phase functional and under suitable conditions on $f_1$ and $g$, it has a non trivial minimizer, say  $U^+ \geq 0$. Then, from  above lemma, we get
\[
J(\min\{\vp, U^+\}) + J(\max\{\vp, U^+\})  \leq J(\vp) + J( U^+).\] Since  $\max\{\vp, U^+\} = \max\{\vp^+, U^+\}$ and $U^+$ is a minimizer for the one-phase functional $J^1_{\{f_1, g\}}$, then it follows that
\be
J(\min\{\vp, U^+\}) \leq J(\vp) ~~~~{\rm for~any~~~~} \vp \in H^1(\R^N).\ee
In particular, this shows that the minimizer of the functional $J_{\{f_1, f_2,  g \}}$, if it exists, can be assumed to have support inside the union of supports of the minimizers for the corresponding one phase functionals, $J^1_{\{f_1, g\}}$ and $J^1_{\{f_2, g\}}$ .

\begin{corollary}\label{cor1} Suppose that the functional $J_{\{f_1,f_2,g\}}$ has a minimizer $u \in H^1( \R^N )$, with nonempty supports $ \Om^\pm := \{ x \in \R^N : \pm u(x) > 0 \}$. Let $U^1\geq 0$, $u^1 \leq 0$ be functions such that
 \bea
J_{\{f_1,g\}}^1 (U^1) & := &   \inf\limits_{ \{\vp \in H^1(\R^N) : \vp \geq 0 \}} J_{\{f_1,g\}}^1 (\vp) \\
J_{\{f_2,g\}}^1 (u^1) & := &   \inf\limits_{ \{\vp \in H^1(\R^N) : \vp \leq 0 \}} J_{\{f_2,g\}}^1 (-\vp), \eea
 where
\be
J^1_{\{f,g\}} =   \int\limits_{\R^N} \left( |\nabla u|^2 - 2f  u  + g^2 \chi_{\{u > 0\}} \right)  \, dx
\ee is the one phase functional.
 Then, we conclude that
\be
\Om^+ \subseteq \{ U^1 > 0 \},~~~~~ \Om^- \subseteq \{ u^1 < 0 \}.\ee
\end{corollary}
 In fact, for any $\vp \in H^1(\R^N)$, we can write
\be\label{00}
J(\vp) = J^1_{\{f_1,g\}}(\vp^+) +  J^1_{\{f_2,g\}}(-\vp^-) \geq \min\limits_{ v \geq 0} J^1_{\{f_1,g\}}(v) + \min\limits_{ v \leq 0} J^1_{\{f_2,g\}}(v).
\ee
Thus, as long as there are conditions on $f_1$, $f_2$ and $g$ which guarantee existence 
of one phase minimizers with compact support, we can minimize the functional $J_{\{f_1, f_2, g\}}$ on $H^1(\R^N)$. 
In particular, we recall Theorem 1.4 of \cite{gs}  
which give conditions for existence of minimizer for  the one-phase functional $J^1_{\{f,g\}}$, viz.,
 \be\left.\begin{array}{lll}
 (A 1)&& f, g \in L^\infty(\R^N) \\
 (A 2)&& {\rm supp\/} ~f^+~{\rm~ is ~compact}\\
 (A 3)&& g \geq 0\\
 (A 4)&& {\rm at~ least~one~of~} f \leq - c_1< 0{~\rm  \ or  \ \/}~~~ g \geq  c_0 > 0~\\
  &&~~~~~~~{\rm hold~outside~a~compact~set~for ~ some~positive~constants\/~}c_0,~c_1.
 \end{array} \right\}\ee
 The following theorem gives existence for two phase functional.
 \begin{theorem}\label{exist}(Existence) 
Consider functions $f_1$, $f_2$, $g$  satisfying  (A1) and such that
\bea
&& {\rm supp\/} ~f_1^+~{\rm and~}  {\rm supp\/} ~f_2^-~{\rm~ is ~compact}\\
&&   f_1 \leq -c_0 < 0~ {\rm and \/}~~~ g^- \geq c_1 > 0~{\rm hold~outside~a~compact~set\/} \\
&& {\rm or,~} f_2 \leq  -\tilde c_0 < 0~ {\rm and \/}~~~ g^+ \geq \tilde c_1> 0~{\rm hold~outside~a~compact~set\/},
\eea
where $c_0$, $\tilde c_0$, $c_1$ and $\tilde c_1$ are positive constants. Then, there exists a minimizer for the functional $J_{\{f_1,f_2,g\}}$ in $H^1(\R^N)$. 
 \end{theorem}
\begin{proof}

  Since  $(f_1, g^+)$ and $(-f_2, - g^-)$ both satisfy the conditions $(A1)-(A4)$, then we get existence of minimizers for the one phase functionals  $J^1_{f_1, g}$ and $J^1_{f_2, g}$.
Thus, minimizing $J_{\{f_1, f_2, g\}}$ over the set 
\[
W:=\{ u \in H^1(\R^N) : u^1 \leq u \leq U^1 \}\]
and repeating the proof of Proposition 2.1 of \cite{bjh},  we obtain a minimizer for the two phase functional $J_{\{f_1, f_2, g\}}$. Here we note that 
$J_{\{f_1, f_2, g\}} \geq J_{\{f_1, f_2, g = 0 \}}$. 

\end{proof}

Theorem \ref{exist} will be used to prove the existence of "two phase quadrature surface" in Section 7.
This is the case  when each $f_i$ in \eqref{j2} is replaced by $\mu_i \ast \psi$, where $\mu_i\ast \psi $ is a mollified version of  a  positive Radon measure $\mu_i$ with compact support. We will restrict ourselves to the case when measures  $\mu_i$  are "sufficiently concentrated" as defined in \cite{gs}, which we refer to as measures satisfying Sakai's concentration condition defined as follows.

 \begin{definition}[Sakai's  concentration condition]\label{sakai}
 	We say that the Radon measure $\mu$ satisfies \emph{Sakai's concentration condition} if for every $x\in supp(\mu)$
 	 	\[
 	\underset{r\to 0^+}{\;\limsup}\;{\frac{r\cdot\mu(B_r(x))}{|B_r|}}> \frac{N6^Nc}{3},
 	\]
 	where $c>0$ is a fixed constant such that $0\leq g(x)\leq c$. Here $g$ is the function given in  \eqref{j2}.
 \end{definition}


\section{Free boundary condition for the  Minimizer}
\setcounter{equation}{0}

The {\em free boundary}  of $u$, denoted $\Gamma_u=\Gamma$, is defined as 
$\Gamma_u=\partial \Omega^+ \cup \partial \Omega^- $ where
$$
\Omega^\pm=\Omega^\pm_u=\{x\in \R^N :\; \pm u(x)>0\}.
$$
A point $ x \in \Gamma_u$ is said to be a {\em one phase free boundary point\/}  if there exists  $r> 0$ such that 
\[ 
\overline{\Omega^+}\cap\overline{\Omega^-}\cap B_r(x)=\emptyset \]
and it is said to be a {\em two phase free boundary point\/} if for all $r > 0$,  
\[\overline{\Omega}^+\cap \overline{\Omega}^-\cap B_r(x)\ne\emptyset. \] 
 The {\em set of one phase free boundary points} of $u$ is  defined as 
$$
\Gamma'=\{x\in \Gamma_u : \textrm{ there exist an }r>0\textrm{ such that }\overline{\Omega^+}\cap\overline{\Omega^-}\cap B_r(x)=\emptyset\}.
$$
while {\em the set of  two phase points}, denoted $\Gamma''_u=\Gamma''$ is 
$$
\Gamma''=\{x\in \Gamma;\; \textrm{ for all }r>0\textrm{ we have }\overline{\Omega^+}\cap \overline{\Omega^-}\cap B_r(x)\ne\emptyset \}.
$$
Finally, {\em the set of branch points\/} $\Gamma^*=\Gamma^*_u$ is the intersection of $\overline{\Gamma'}$
and $\overline{\Gamma''}$;
$$
\Gamma^*=\overline{\Gamma'}\cap \overline{\Gamma''}.
$$
The free boundary $\Gamma_u$ of a solution $u$ to a two-phase problem 
 can thus be   decomposed as 
\be\Gamma_u = \Gamma' \cup \Gamma'' \cup  \Gamma^*.
\ee

Here we show that under suitable conditions, a  local  minimizer $u$ of the functional $ J
:= J_{\{f,g\}}$ satisfies 
\be 
\Delta u + f {\ml } \mid_{\Om} = g{\mh}^{N-1} \mid_{\pa \Om}\quad {\rm in }~  \Om:= {\rm supp} (u).\ee
 Observe that, for any $ \vp \in C^2_0(\Om)$,  $\chi_{\{ u + t \vp \neq 0 \}} = \chi_{\{ u \neq 0 \}} $ and  we have
\Bea
 && J(u + t \vp) - J(u) \\
& = &   \int\limits_{\R^N} \left\{ |\nabla (u + t \vp) |^2 - 2f_1 (u + t \vp)^+ + 2f_2 (u + t\vp)^- + g^2 \chi_{\{u + t \vp \neq 0\}} \right\}  \, dx \\
  &&~~~~~~~- \int\limits_{\R^N} \left\{ |\nabla u|^2 - 2f_1 u^+ + 2f_2 u^- + g^2 \chi_{\{u \neq 0\}} \right\}  \, dx \\
& =  & \int\limits_{\R^N}   \left\{ 2 t \nabla u \cdot \nabla \vp + t^2 |\nabla \vp|^2 - 2f_1 [(u + t \vp)^+ - u^+]  + 2[f_2 (u + t\vp)^-  - u^-]  \right\}  \, dx
\Eea
It follows that $ \vp \mapsto \int\limits_{\R^N} \nabla u \cdot \nabla \vp  - 2 f_1 \chi_{\{u > 0\}} \vp  + 2 f_2  \chi_{\{u <  0\}}\vp $ is the Euler-Lagrange equation for the functional and after integration by parts we have that $u$ satisfies
\be
\Delta u =  -f_1 \chi_{\{u > 0\}}  + f_2  \chi_{\{u < 0\}} \quad {\rm in}\quad \Om .
\ee

One may now show in a standard way that the minimizer satisfies the Bernoulli boundary condition in a weak sense, and in the strong sense  $H^{n-1}$-almost everywhere  on the free boundary. 
Here weak sense  refers to 
 \be\label{weak}
\lim_{\var \searrow 0} \int\limits_{\pa \{\pm  u > \var \} \cap B_r(z^1)  } \left( |\nabla u|^2 - g^2 \right) \Theta \cdot \nu \,d \mh^{N-1}  = 0 
\ee
where $z^1 \in \Gamma_u'$ with  $B(z^1, r) \cap \overline{\Om^-} = \emptyset$ 
(respectively, $B(z^1, r) \cap \overline{\Om^+} = \emptyset$), and
for all {\em vector fields \/} $\Theta \in C^0( B_r(z^1), \R^N )$. Here $\nu$ denotes the outward normal vector to the boundary of the sets. Analysis of the free boundary in neighborhood of the branch points is very technical and relies on the results proved in the paper \cite{ASW}. The following theorem summarizes the regularity properties of free boundary:

 \begin{theorem}(\cite{ASW})
Let  $u$ be a minimizer of $J_{\{f_1,f_2, g\}}$. Then following holds. 
\begin{itemize}
\item[i)] The weak free boundary condition \eqref{weak} holds for minimizers.
\item[ii)] For any point $z \in \Gamma'' \setminus \Gamma^\star $ (two-phase and non-branch points) we have $\Delta u = 0$ in $B_r(z)$, provided  $B_r(z) \cap \Gamma' = \emptyset$.
\item[iii)]  The free boundary has finite $(N-1)$-dimensional Hausdorff measure.
\item[iv)] Close to the branch points, the free boundary consists of two $C^{1, \alpha}$ graphs in a universal neighborhood of the branch point.
\end{itemize}
\end{theorem}
The proof of (i) is straightforward and  similar to that of Theorem 2.4 in \cite{ACF}, by use of domain variation. Note that our model, as formulated,  does  not require the condition ${\rm meas} \{ u = 0 \} = 0$. This is because in our situation,  the function $\la(u) $ of \cite{ACF} is 
\be
 \la(u)  =  g^+ \chi_{\Omega^+} + g^- \chi_{\Omega^-}
  \ee 
 and hence $\la(0) = 0 $.  From Theorem 7.1 and Remark 7.1 in \cite{ACF}, it follows that the set of  one phase boundary points $\Gamma_u'$ has finite $ (N-1)$-Hausdorff measure.  Furthermore, due to the choice $g_1=g_2$ we have made here, the two phase free boundaries are level surfaces of harmonic functions.
 
 Proof of iii) also uses non-degeneracy of both phases.\footnote{This is in general not true for the functional in \cite{ACF}. For our case we have the advantage of $\lambda (0)=0 $. } For interested reader we refer to \cite{ASW} Section 3, for further local measure theoretic properties of the free boundary.
 
 The proof of iv) is a deep result, using chains of technical arguments. The core idea is that due to non-degeneracy of both phases, a blow up (scaling of the type $u_r(x):=u(rx+ z) /r$ at any branch point $z$) 
 leads to a global two-phase solution, which is classified and shown to be a two-plane solutions (i.e. a broken linear function  $L(x):=a_+ x_+ - a_-x_-$). One may then reiterate the blow-up argument, but this time by linearization technique, i.e. considering the limits of $(u_r - L(rx) )/r$. One proves that these limits exists and will solve a so-called thin-obstacle problem, which in turn is well-studied. Enough information and  knowledge about the regularity theory of  their free boundary is available in literature.  From here on, one may then show that our model is a (close enough) perturbation of the limiting problem and hence we can derive regularity of the free boundary for our original problem.
 
 \section{Qualitative properties}
 In this section we discuss  some qualitative properties  of minimizers of (\ref{j2}), that has already been established for the one-phase case. 
 The complications, with  the two-phase case makes similar properties much more hard to show.   Here we apply the  moving plane method to obtain convexity or monotonicity of the level sets of minimizers. To this aim, for a fixed unit vector  $n\in \R^N,$ and for $t\in\R $ we set
 \[
 T_t=\{x\cdot n=t\}, \;\; T^-_t=\{x\cdot n<t\}, \;\;\mbox{and}\;\;T^+_t=\{x\cdot n>t\}.
 \]
 For $x\in \R^N$ let $x^t$ denote the reflection of $x$  with respect to $T_t$. We also set $\varphi^t(x)\equiv\varphi(x^t)$,  for a function $\varphi$ and if  $\Omega\subset\R^N $ we define
 \[
 \Omega_t=\Omega\cap T^+_t \;\; \mbox{and} \;\; \tilde{\Omega}_t=\{ x^t :\; x\in\Omega_t \}.
 \]
 \begin{theorem}\label{geometr_prop_m=2}
 	Let $f_i(x)$, $i=1,2$ satisfy conditions $(A1)-(A4),$ and assume that for some unit vector $n\in \R^N$ and some $t_0\in\R $ we have
 	\be
 	f_i(x)\leq f^t_i(x),\;\; g(x)\geq g^t(x)\;\;\mbox{in}\;\; T^+_t,\quad 
 	\mbox{for~ all}~ t\geq t_0.  
 	\ee
 	Then for a smallest minimizer $u \in H^1(\R^N)$ of the functional ${J}_{\{f_1,f_2,g\}}$, we have
 	\begin{align*}
 	u<{u}^t\;\;\mbox{in}\;\; \Omega_t\;\;&\mbox{for all}\;\; t\geq t_0, \\
 	\tilde{\Omega}_t\subset\Omega\;\;&\mbox{for all}\;\; t\geq t_0\\
 	\end{align*}
 	where $ \Omega =  {\rm int}( \overline{\Om^+\cup\Om^-})$, with $\Om^\pm :=\{ x \in \R^N: \pm u(x) > 0\}$. 	\end{theorem}
 \begin{proof} For $t \geq t_0$, set
 	\[
 	v^t=
 	\begin{cases}
 	\min(u,{u}^t),\;\;\mbox{in}\;\; T^+_t,\\
 	\max(u,{u}^t)\;\;\mbox{in}\;\; T^-_t.
 	\end{cases}
 	\]
 	Let
 	\[
 	L(\varphi)=\int_{T^+_t} \left(|\nabla \varphi|^2-2f_1\varphi^++2f_2\varphi^-+g^2\chi_{\{\varphi\neq 0\}}\right)dx,
 	\]
 	and
 	\[
 	L_t(\varphi)=\int_{T^+_t} \left(|\nabla \varphi|^2-2f^t_1\varphi^++2f^t_2\varphi^-+(g^t)^2\chi_{\{\varphi\neq 0\}}\right)dx.
 	\]
 	According to Lemma \ref{l1} we obtain
 	\begin{align*}
 	{J}_{\{f_1,f_2,g\}}(v^t)&=L(\min(u,{u}^t))+L_t(\max(u,{u}^t))\\
 	&\leq L(u)+L_t({u}^t)={J}_{\{f_1,f_2,g\}}(u),
 	\end{align*}
 	for all $t\geq t_0.$ Thus,$${J}_{\{f_1,f_2,g\}}(v^t)={J}_{\{f_1,f_2,g\}}(u).$$ For all $t\geq t_0$ so large that $\Omega\subset T^-_t$ we  apparently have $v^t=u,$ i.e. $u\leq u^t$ in $T^+_t.$ This yields $\tilde{\Omega}_t\subset\Omega$ for all $ t\geq t_0.$  Since $u$ is a smallest minimizer, then we conclude $u\le v^t,$ which completes the proof.
 	
 \end{proof}
 
 \begin{corollary}\label{c1}
 	Let  ${u}\in H^1(\R^N),$ and $f_i(x), i=1,2$ be as in Theorem \ref{geometr_prop_m=2}. If we assume  $f_i(x), i=1,2$ are symmetric in $T_{t_0},$ then minimizer $u$ is symmetric in $T_{t_0}.$
 \end{corollary}
 
 \begin{corollary}\label{c2} 
 	Let $u$ be a solution of 
 	\be\label{locmin}
 	\begin{cases}
 		\Delta u = -f_1 \chi_{\{u>0\}} + f_2 \chi_{\{u<0\}}\quad &{\rm in} ~\Om, \\
 		u = 0, |\nabla  u| =  g \quad &{\rm on } ~\pa \Om,
 	\end{cases}
 	\ee
 	with  
 	\be \Omega =  {\rm int}( \overline{\Om^+\cup\Om^-}),
 	~~~\Om^\pm :=\{ x \in \R^N: \pm u(x) > 0\} .\ee
 	Suppose that  $\mu = c_+\delta_{z^+} - c_-\delta_{z^-}$ is Dirac measure and $g^2 \equiv constant> 0$. Then the solution of (\ref{locmin}) is symmetric with respect to the line joining the points $z^+$ and $z^-$.
 \end{corollary}
 \begin{proof} The Corollary \ref{c2} follows from Theorem \ref{geometr_prop_m=2} by choosing $f_1 = c_+ \rho$ and  $f_2 = c_- \rho$ where  $\rho$ is a radially symmetric  mollification of Dirac measure. It follows that $\Omega $ is  has rotational symmetry with respect to the axis $L:=z^- + t (z^+ - z^-)$ ($t \in \R$), and $|u| (x) \geq  |u|(y)$ for $(x-y)$ orthogonal to $L$ with   $d(x, L) < d(y, L)$. 
 	
 \end{proof}
 
 \begin{remark}
 	Observe that  symmetry of  $\Om^+$ and/or $\Om^-$ will depend on the weights $c^+$ and $c^-$. 
 	In particular, if say  $\mu^\pm$ is sufficiently concentrated around the point $z^\pm$ so that $\overline{ \Om^+} \cap \overline{ \Om^-} = \emptyset $,  then  $\Om^\pm$ is a ball and $u^\pm$ is bounded,  radially symmetric. 	
 \end{remark}


\section{Quadrature Identities}\label{sec-qi}

\subsection{Two-phase Quadrature surfaces}
In this section we  discuss the concept  of two-phase Quadrature surfaces.
The one phase problem has been well studied in the literature and we refer the reader to the following works \cite{gs,henrot,onodera,h2,h1}. 

As for two-phase quadrature domains case (see \cite{bjh}), here again the key point is that the measures have to be concentrated enough and also in balance. Indeed, if  the measures $\mu_1$, say, has a very high density on its support, but not the other, then the support of the corresponding $u_1$ will have the possibility of  covering  the support of $\mu_2$.  This naturally makes it impossible to find a two-phase QS for our measures. 
Finding right conditions for this balance is a question to be answered in the future. Here we will illustrate this for measures that satisfy Sakai's concentration  condition. 

Let  $f_i(x)=\mu_i-\lambda_i(x)$ for $i=1,2$  where $\lambda_i\geq 0$ are $L^\infty$ functions, and $\mu_i$ are positive bounded Radon measures.  The main difficulty is to provide  conditions which lead to the
 existence of solution of the system \eqref{twophase} with property $supp(\mu_1)\subset \{{u}>0\},$ and $supp(\mu_2)\subset \{{u}<0\}.$ The latter property  implies the following conditions $$
 \mu_1\equiv\mu_1\chi_{\{{u}>0\}}\;\;\mbox{and}\;\; \mu_2\equiv\mu_2\chi_{\{{u}<0\}}.
 $$  Then the  system \eqref{twophase}, for $f_i=\mu_i-\lambda_i,$  can be rewritten as follows:
 \begin{equation}\label{system2}
 	\begin{cases}
 	\Delta u=(\lambda_1\chi_{\{u>0\}} -\lambda_2\chi_{\{u<0\}})-(\mu_1-\mu_2)&\;\;\mbox{in}\;\;
 	\Omega,\\
 	u=0,\;|\nabla u|=g & \;\;\mbox{on}\;\;\partial\Omega,
 	\end{cases}
 \end{equation} 
 where  \be \Omega =  {\rm int}( \overline{\Om^+\cup\Om^-}),
 ~~~\Om^\pm :=\{ x \in \R^N: \pm u(x) > 0\}.\ee

 For quadrature surfaces we need to take $\lambda_1=\lambda_2=0.$ Then for arbitrary  $h \in HL^1 (\Omega^+ \cup \Omega^-),$ we write Greens second identity:
  \begin{equation}\label{green_idnt}
  \int_{\Omega^+ \cup \Omega^- }(h\Delta u-u\Delta h)dx=\int_{\partial(\Omega^+ \cup \Omega^-)}\left(h \frac{\partial u}{\partial \nu}-u\frac{\partial h}{\partial \nu}\right)d\sigma_x.
\end{equation}
 
 Apparently equation \eqref{green_idnt} leads to 
 $$
 \int_{\Omega^+ \cup \Omega^- }h\Delta udx =\int_{\partial(\Omega^+ \cup \Omega^-)} h \frac{\partial u}{\partial \nu}d\sigma_x.
 $$
 Thus, one can formally write (leaving the verification to the reader) 
 \begin{align*}
\int_{\Omega^+ \cup \Omega^-}h(-(\mu_1-\mu_2))dx&=\int_{\partial(\Omega^+ \cup \Omega^-)} h(x) (\nabla u\cdot \nu)d\sigma_x\\&=-\int_{\partial\Omega^+} gh(x)d\sigma_x+\int_{\partial\Omega^-} gh(x)d\sigma_x,
 \end{align*}
 which finally gives 
 \[
\int h(d\mu_1-d\mu_2)=\int_{\partial\Omega^+} gh(x)d\sigma_x-\int_{\partial\Omega^-} gh(x)d\sigma_x.
 \]
 It is easy to see that the standard  mollifier technique (see \cite{gs}) will also work in this case, and we may replace the measures with smooth functions, with support close to the support of measures.

 \begin{definition}[Two-phase Quadrature surface]
 	Suppose we are given two bounded positive measures $\mu_1,\mu_2$ and disjoint domains $\Omega_1,\Omega_2$ such that $supp(\mu_i)\subset \Omega_i.$ If for  harmonic functions $h\in HL^1(\overline{\Omega_1\cup\Omega_2})$ the following QI holds
 	\begin{equation}\label{Q-identity}	
 		\int_{\partial \Omega_1} g h (x) \  d\sigma_x - \int_{\partial \Omega_2} g h (x) \  d\sigma_x= \int h d(\mu_1-\mu_2) \ ,
 	\end{equation}
 	then we call $\partial\Omega_1\cup\partial\Omega_2$ a  {\it  Two-phase QS} with respect to the measure $\{\mu_i\}_{i=1}^2,$ and a given  smooth positive function $g.$
 	
 	If we reduce the  test class $h$ to be subharmonic in $\Omega_1$ and super-harmonic in $\Omega_2$ (due to negative sign in front of the integral) then  the equality  in \eqref{Q-identity} is replaced with an inequality $(\geq) $.
 	 \end{definition}
 	Observe also if we take $\mu_2=0,$   then  $\overline{\Omega}_2=\emptyset,$ and we get the definition of a one-phase quadrature surface
 	
 	\begin{equation}\label{1-Q-identity}
 		\int_{\partial\Omega_1}gh(x)\; d\sigma_x  = \int h d\mu_1.
 	\end{equation}

\begin{theorem}\label{existenceQS}
	Let $\mu_1$ and $\mu_2$ be  given Radon measures with compact supports, that satisfy Sakai's concentration condition as in Definition \ref{sakai}. Suppose that for each $\mu_i$ the corresponding  one-phase quadrature surface
	$\partial Q_i$ (see \eqref{1-Q-identity})  is such that
	\begin{equation}\label{disjointQS}
		\overline{Q}_1\cap supp(\mu_2)=\emptyset,\;\; \mbox{and}\;\;\overline{Q}_2\cap supp(\mu_1)=\emptyset .
	\end{equation}
	Then, we have a solution to our two-phase free boundary problem \eqref{system2} along with  $supp(\mu_1)\subset supp(\{u>0\}),$ and  $supp(\mu_2)\subset supp(\{u<0\}).$
\end{theorem}

\begin{proof}
	We  consider  mollifiers  $\mu_i \ast \psi $, of the measures   $\mu_i$  ($i=1,2$) and minimize the functional  \eqref{j2} for $f_i = \mu_i \ast \psi$. Since $supp(\mu_i \ast \psi ) $ is a subset of a $\epsilon$-neighborhood of $supp(\mu_i)$, it suffices (by taking $\epsilon $ arbitrary small) to show the theorem for smooth $\mu_i$. We thus from now on assume $\mu_i$ is smooth enough such that a minimizer 
	${u}$ belongs to  $ H^1(\R^N).$
	
	Next, if we prove that $supp(\mu_1)\subset supp(\{u>0\}),$ and  $supp(\mu_2)\subset supp(\{u<0\}),$ 
	then obviously ${u}$ will solve the two-phase free boundary problem \eqref{system2}. To this end, first observe that due to Sakai condition for QS, we have the following embedding $supp(\mu_i)\subset Q_i,$ for $i=1,2.$ We argue by contradiction. Assume that $supp(\mu_{1})\setminus supp(\{u>0\})\neq\emptyset.$ The case for the measure $\mu_2$ can be done in a similar way. Then, according to condition \eqref{disjointQS} and Corollary \ref{cor1}, there exists a point  $z_0\in supp(\mu_1)\setminus supp(\{u>0\}),$ such that $dist(z_0,\tilde{\Omega})>0,$ where $\tilde{\Omega}=supp(\{u>0\})\cup supp(\{u<0\}).$ Thus, one can easily take a ball $B_R(z_0)$ such that $\overline{B_R(z_0)}\cap\tilde{\Omega}=\emptyset.$
	
	Let the constant $M>0$ is taken such that $\max\left(\underset{i=1,2}{\sup}\;|\mu_i|,\;\cfrac{Nl_0}{R}\right)< M<+\infty,$ where $g(x)\geq l_0>0.$  
	We set $$r=R\cdot\left(\frac{Nl_0}{MR}\right)^{1/N},$$ and consider the following measure, which  satisfies Sakai's condition
	\[
	\nu_{1}\equiv\mu_{1}\cdot\chi_{B_{r}(z_0)}.
	\]
	Define
	\[
	L_{M,l_0}(U)=\int\limits_{\R^N}\left(|\nabla U|^2 - 2M\cdot\chi_{B_r(z_0)}U+l_0^2\chi_{\{U>0\}} \right)dx.
	\]
	According to Lemma $1.2$ in \cite{gs} every minimizer of $L_{M,l_0}(U)$ over the set $ \{\vp \in H^1(\R^N) : \vp \geq 0 \}$ is radially symmetric, radially non-increasing and vanishes outside a compact set. Moreover,  the largest minimizer support (in our case $c>0,\; b=0$ and $R_1=0$)  is a ball centered at $z_0$ and with radius $\sigma =\left(\frac{r^NM}{Nl_0}\right)^{1/{N-1}}=R$ (see Example $1.5$ and the proof of Theorem $1.4$ in \cite{gs}).  The proof of this result relies on the so-called symmetric decreasing rearrangement technique, and we refer for its background to the book \cite{mossino}. 
	
    Let $v_1$  be a minimizer to the functional
	\[
	J^1_{\{\nu_{1},g\}}(U)=\int\limits_{\R^N}\left(|\nabla U|^2 - 2\nu_1 U+g^2\chi_{\{U>0\}}\right) dx,
	\]
	over the set  $ \{\vp \in H^1(\R^N) : \vp \geq 0 \}.$
   Then using the same arguments as in the proof of  Theorem $1.4$ in \cite{gs}, one can easily conclude that   $supp(v_1)\subset \overline{B_R(z_0)},$ and therefore $supp(v_1)\cap\tilde{\Omega}=\emptyset.$
	 Apparently $supp(\nu_1)\subset\{v_1>0\},$ which implies that $J^1_{\nu_1,g}(v_1)<0.$
	Now, simple computation gives that
	$$
	J_{\{\mu_1, \mu_2,g\}}(u+v_1)\leq J_{\{\mu_1, \mu_2,g\}}(u)+J^1_{\{\nu_1, g\}}(v_1)<J_{\{\mu_1, \mu_2,g\}}(u).
	$$
	This contradicts the minimality of $u.$ Thus $ supp(\mu_{1})\setminus supp(\{u>0\})=\emptyset,$ and this implies  $supp(\mu_{1})\subset supp(\{u>0\}).$ This completes the proof of Theorem
\end{proof}


\subsection{Examples of two-phase QS}
\setcounter{equation}{0}
 Due to Theorem \ref{exist},  minimizers for $J_{\{f_1, f_2, g\}}$ always exist in the following  special cases:
\begin{itemize}
	\item[a)]  $f_1 \equiv 0$ and  both  $g^+ >0$,  $   f_2 \leq c_2 < 0$ outside a compact set or  $f_2 \equiv 0$ and  both  $g^- >0$,  $   f_1 \leq c_1 < 0$ outside a compact set.
	\item[b)] A simpler two phase functional where it is assumed that either $g \geq 0$ or $g \leq 0$. 
\end{itemize}
These conditions, however, do not say anything about whether  the quadrature identity \eqref{qi} is admitted by the domain and the measure.  Here we discuss  simple examples of two phase QS, some of which are generated using one phase QS and symmetry arguments.

\noindent {\bf Example 1 (Plane Symmetric QS): \/} 
Let $(u,\Omega)$ be a one phase QS for a measure $f$, and $g$ as appearing in the functional \eqref{j1}, see also equation \eqref{od}. Consider  further  a hyperplane $T$,  not intersecting $\Om $,  and  an 
odd reflection of the solution $u$  with respect to a plane $T$.
This  will trivially give a two phase quadrature surface 
\be\left\{
\begin{array}{cll}
	- \Delta u & =   f \chi_{\{u > 0\}}  + \tilde f  \chi_{\{u < 0\}} \quad & {\rm in}\quad \Om \cup  \tilde \Om ,\\
	u & =  0,\quad & {\rm on }  \quad \pa \Om \cup \tilde \Om ,\\
	\frac{\pa u}{\pa \nu} &=  1  \quad & {\rm on } \quad \pa \Om,\\
	\frac{\pa u}{\pa \nu} & =-1 \quad  & {\rm on } \quad \pa ( \tilde \Om ),
\end{array} \right.
\ee
where $\tilde \Omega  $, $\tilde f$ denotes the reflection  of $\Omega$, respectively  $f$, in the plane  $T$.

A less trivial example can be constructed as follows: 
Let $g$, and $f$  be as before for the one phase QS. Let further $T^+:=\{ x : x_1 > 0\}$ and minimize the one-phase functional in the set $T^+$, with zero Dirichlet boundary values on $\partial T^+$. Suppose further that the support of $u$ reaches all the way to the plane  $\partial T^+$. This will formally solve
\be\label{half-space}\left\{
\begin{array}{ll}
	\Delta u = -f  \qquad &\hbox{in  }\Om = \{ u > 0 \} \subset T^+,\\
	u =0  \qquad &\hbox{in  } \partial T^+ ,\\
	|\nabla  u| = g \qquad &\hbox{on } \pa \Om .
\end{array} \right.
\ee
Then odd reflection of $u$ with respect to the plane $T$ gives a  quadrature surface symmetric about the plane $T$. 

\noindent {\bf Example 2 (Spherically  symmetric  QS): \/} 
A different example  would be  an annular two phase quadrature surface. That is,  a quadrature surface $\Gamma = \pa \Om= \pa \Om_1 \cup \pa \Om_2$ such that $ \Om= \Om_2 \setminus \Om_1$ is an annular domain with inner boundary $\partial \Om_1$ and outer boundary $\partial \Om_2$ with $g^+ = c_1 >0 $ on $ \pa \Om_1$ and $g^- = c_2 < 0 $ on $ \pa \Om_2$ (or vice versa). For a construction of a spherical annular two phase quadrature surface, we proceed as follows.  Consider a   uniformly distributed (and sufficiently large)\footnote{It suffices to take $\mu = 3d\sigma$, where $d\sigma$ is the surface measure.}  measure  $ \mu$ on  the sphere $S_2 :|x|=2$ (or defined in a $\var$- tubular neighborhood of $S_2$), and  solve the one phase  free boundary problem
\be\label{ann}
\begin{cases}
	- \Delta  u  =  \mu & \;\;{\rm  in}\;\;  B_R \setminus B_1  \\
	u  = 0 &  ~~~~~{\rm  on}~~~~~  |x|=1 ~{\rm and~} |x| =  R\\
	|\nabla u |  = 1 &  ~~~~~{\rm  on}~~~~~   |x|=R. 
\end{cases}
\ee
It is not hard to calculate explicitly what $R > 1 $ is, but we surely know that there is at least one such $R$. 

Now a two phase solution  can be obtained by extending $u$ by  an odd Kelvin inversion  of  $u$ with respect to  the sphere $|x|=1$. 

Then the extended function, which for simplicity  is labeled  $ u $, satisfies  
\be
\begin{cases}
	- \Delta  u  = \nu &  \;\;{\rm  in}\;\; 1/R < |x| < R ,\\
	u  =  0 &  \;\;{\rm  on}\;\;|x|=1/R \;\;  |x|=1, \,  {\rm and}\;\; |x| =  R,\\
	|\nabla u |  = R^{-n-2} &  \;\;{\rm  on}\;\;|x| = 1/R   ,  \\
		|\nabla u |  = 1 &  \;\;{\rm  on}\;\;|x| = R  ,  \\
\end{cases}
\ee
Here $\nu = \mu - \tilde \mu$, where $\tilde \mu$ is the even Kelvin reflection of $\mu$ in $|x|=1$.
Then the boundary of the  new domain is  the required quadrature surface. 

		\noindent {\bf Example 3 (Non-trivial two-phase QS): \/} 
		The above examples illustrates that it is not trivial to give  explicit examples of    QS that do not carry information from the one-phase problem. Here we shall  give one slightly more complicated example, which again is constructed by reflection of one-phase QS. Nevertheless, the reflection is  more elaborate than standard ones, and it is defined  through the so-called Schwarz function \cite{Davis}, which is defined as follows:\\
		Let $\Gamma$ be any  analytic curve dividing the complex plane $\mathbb C$ into at least two components.
		Denote by $\Gamma^+$ one of these components. 
		The Schwarz function, which we denote by $S(z)$ for $z=x + {\bf i}y$, is a function which is analytic in a neighborhood, say $N$, of $\Gamma$ satisfying $S(z) = \bar z$ on $\Gamma$. Next fix a point  $z^0 \in \Gamma$, and  suppose that $B_r (z) \subset N$ so that $S(z)$ is analytic in $B_r (z)$.
		We assume without loss of generality that $r=1$, otherwise we just scale the curve with $r$.

		Choosing  $g=1$ and $\mu$ a "smooth" measure  with support in the open set $\Gamma^+ \cap B_1 (z) $, let  $(u,\Omega) $ be a minimizer of the functional $J_{\{\mu, g\}}$  with zero Dirichlet data on $ \Gamma \cap B_1 (z) $.  We may further assume that  $\mu$ satisfies conditions so that there exists a  solution to the one-phase equation \eqref{od} in $\Gamma^+ \cap B_1(z)$, cf. equation 	\eqref{half-space}. 	Next, we  invoke the anti-conformal reflection  ${\mathcal R}_{\Gamma} $ associated to $\Gamma$,  and defined as ${\mathcal R}_{\Gamma}  (z)= \overline{ S(z) }$ (see \cite[Chapter 6]{Davis}).
		
		For $z\in \Gamma^- $ we define  $u(z)=-u({\mathcal R}_\Gamma z)$, and in this way we extend $u$ across $\Gamma$ as a solution to our problem with negative measure $-\mu({\mathcal R}_\Gamma z )$.	This creates an example of a  two-phase free boundary  for our problem.

\subsection{Solutions with unbounded support} 

There are not many trivial examples of  two-phase QS with unbounded support,  however, there exist a few. The most simple example is naturally when we take $\mu_1$ to be Dirac mass at origin
and $\mu_2 \equiv 0$. Then the appropriate sphere is both one- and two-phase QS, that can easily be verified, using
integration by parts. Continuing on this path, if we 
 assume both the  measures to be identically zero,  then for  $g_1 = g_2=  \mbox{  constant}$ one can show that an  appropriate  linear function is a solution to our problem.  
Quadrature surfaces, with  unbounded supports (for both phases)  can be constructed  from bounded ones, by a simple procedure. Indeed, if  we already have a QS, for a measure $\mu$, we may consider minimizing the corresponding functional in $B_R$, where we now put Dirichlet data on $|x|=R$,  that corresponds to $h(x)=(x - x^0)\cdot {\bf a} $ for some $x_0$ and vector ${\bf a}$, such that $supp (\mu) \subset \{ h > 0 \}$.  Any (global) minimizer $u_R$ to this problem will have the property that  its support stretches all the way up to sphere $|x|=R$, due to the boundary values.
Such domains give rise to (partial) QS, which amounts to being QS for the class of harmonic functions on  $supp (u_R)$, that vanishes on the sphere. By letting $R$ tend to infinity, along with using barrier arguments for control of the (linear growth) one can show that there is a limit (at least for a subsequence of $R$) which satisfies a quadrature identity. The heuristic argument here can be made easily rigorous by some footwork, and is left to interested reader. For quadrature domains 
there are at least two references the authors are aware of \cite{bahrami,Sakai-lecture}.
Similar methods can be applied to a two-phase QS, without much efforts.  

It is interesting to mention that unbounded two-phase QS may behave much more differently than their one-phase counterpart. Indeed, we expect that   two-phase quadrature surfaces, with both phases being unbounded, have to behave like plane solutions at infinity. This can be seen easily if
the QS is smooth, by shrinking the solution through any sequence $u_j = u(R_jx)/R_j$ and obtaining a new unbounded solution $u_\infty$, without any measure (these are called Null QS). One can then classify 
Null QS, which are solutions to $\Delta u_\infty = 0$ outside the zero set of $u_\infty$, and have the property that $|\nabla u^+|^2 -  |\nabla u^-|^2 = \hbox{constant}$, where the constant is the limit value of $g_1^2 - g_2^2$ at infinity (see Remark \ref{rem-1}). Since in our case  we have taken $g_1 = g_2$, this implies that we actually obtain a limiting function that is harmonic  with linear growth, and hence a plane. This proves our claim.  For more general values of $g_1, g_2$, one may still prove a similar result, but that would require using strong tools, such as monotonicity formulas, which is outside the scope of this paper.

We close this paragraph by remarking that  bounded two-phase QS are uniformly  bound. This follows from the fact that  two-phase QS are smaller than the union of the corresponding two one-phases, which in turn are uniformly bounded. Again, the details are left to the readers.

\subsection{Uniqueness}

In \cite{h1}, it was shown that $ \pa \Om$ is a quadrature surface with respect to the measure $\mu$ if and only if there is a solution to the Cauchy problem 
\be
    \begin{cases}
	-\Delta u  =    \mu  & \quad {\rm in ~} \Om, \\
	u  =  0 & \quad {\rm on ~} \pa\Om, \\
	\frac{\pa u }{\pa \nu }  = -1  & \quad {\rm on ~} \pa\Om ,
	\end{cases}
\ee
where $ \nu$ is an outward normal. Furthermore, it was proved that if $ \mu = c \delta_x$, $c>0$ and $ \delta_x$ Dirac measure then $\pa \Om$ is a sphere centered at $x$. 

Uniqueness for QS in general fails, unless one has some geometric restriction. This is already known for the 
one-phase problem. Since the functional representing the problem is not convex,  one expects that  local minima as well as stationary points  may give rise to solutions to our free boundary problems.\footnote{It should be remarked that there are other methods such as Perron's smallest super-solution, singular perturbation, and implicit function theory, that have successfully been applied for the existence of one-phase Bernoulli problem, which maybe seen as a QS, when the Dirichlet data is replaced by a source. } For the one phase case there are indeed examples of non-uniqueness for QS worked out by A. Henrot \cite{henrot}.
Therefore, a uniqueness question is even more  complicated in the two-phase case and it seems that the only way to achieve partial results is by imposing strong geometric or other type of restrictions on the solutions, and the data involved.

In one phase problem it is well-known that a QS $(u,\Omega)$ for a single (multiple of) Dirac mass $c_0\delta_{z}$, at the point $z$, is the appropriate sphere $\partial \Omega = \partial B_r(z)$, with $r=r(c_0)$, provided 
$\partial \Omega$ is smooth enough (usually $C^1$ suffices).  The same question for the two-phase problem, for $\mu = c_+\delta_{z^+} - c_-\delta_{z^+}$ seems to be much harder to find an answer to.

\subsection{Null Quadrature surfaces}

In this section we shall let $g\equiv 1$, and discuss the so-called unbounded QS, with zero measures, or so-called null-quadrature surfaces.

A null QS, is a quadrature surface with zero measure (see \cite{karp} for the quadrature domain counterpart). The one phase null-QS refers to a domain $\Omega$ such that one can find a harmonic function $u$  in $\Omega$, with zero Dirichlet data and  $|\nabla u|=g$ on $\partial \Omega$. Obviously $\Omega$ cannot be bounded (due to maximum principle). So one may then wonder about   the  behavior of  $u$  at infinity.

In order to understand the concept of null-QS better, we shall consider it from a potential theoretic point of view, which is more instructive. 
We define, in analogy with null quadrature domains, a null quadrature surface to be the boundary of a domain $\Omega$ such that 
$$
\int_{\partial \Omega} h(x) \ d\sigma_x = 0
$$
for all functions $h$, harmonic  in $\Omega$, and integrable  over $\partial \Omega$.
This is the one-phase case of a null-QS, for $g\equiv 1$.

Let us give a few examples of  one-phase null-QS. The simplest example is the half-space, 
$\Omega =\{a\cdot x >0\}$, ($|a|=1$) with the corresponding function $u=a\cdot x$. A second example is $\Omega= \{a\cdot x >0\} \cup \{a\cdot (x-x^0 <0\}$, where $a\cdot x^0 > 0$.

A third example is the exterior of  any ball $\Omega=\R^n\setminus B_r(x^0)$, where the function $u=b|x-x^0|^{2-n} +c$  (for appropriate $b,c$) solves the
free boundary problem. Naturally cylinders can be built, with exterior of balls as base.

More complicated examples can be given using the construction of H. Alt and L. Caffarelli  \cite{AC}, which is a cone
\[
u(x) := r \max \left \{ \frac{f( \theta)}{f'(\theta_0)}, 0 \right\}
\] using polar coordinates \[
x(r, \vp, \theta) = r (\cos \vp \sin \theta, \sin \vp \sin \theta, \cos \theta)\]
in $\R^3$. The function 
\[
f(\theta) = 2 + \cos \theta \log \left ( \frac{1-\cos \theta}{ 1 + \cos \theta} \right) \]
is a solution of 
\[
(\sin \theta  f')' + 2 \sin \theta f = 0, \quad f'(\frac \pi2) = 0, \]
and $\theta_0 \approx 33.534^\circ$ is the unique zero of $f$ between $0$ and $\frac \pi2$. The function $u$ is harmonic in $\{ u > 0 \}$ and $\pa_\nu u = 1 $ on $\pa \{ u > 0 \} \setminus \{0\}$, i.e.,  the free boundary condition is satisfied everywhere on the surface of the cone, but 
at the origin. At the origin one has $\liminf_{x\to 0}|\nabla u |(x) < 1$. However, since the free boundary is satisfied at every other free boundary point, and that the solution function $u(x)$ behaves linearly at infinity,  one obtains (by simple drill of integrations by parts) 
that $\partial \{u >0 \}$  is a quadrature surface.

Other less regular, and very complicated examples, are the so-called  pseudospheres of John Lewis
\cite{lewis}.\footnote{ John Lewis constructed such objects for Dirac masses, but the same can be done for any measure $\mu$ with high enough concentration, such that the Greens potential 
$G_D^\mu$ of $\mu$ with respect to some  domain $D \subset supp (\mu)$, has the property
 that $|\nabla G^\mu_D| > 1$ on $\partial D$.}
 These objects are much more complicated that fail to be smooth at some points, with gradient of the potential functions being unbounded at some boundary points. Nevertheless, they admit QS  identities, and hence are Quadrature surfaces. 

It is worth mentioning that a recent example of \cite{HHP} (cf. also \cite{traizet})
 solving our  PDE with unbounded support,
has growth that is exponential and that does not qualify as a  QS, in our sense.

The two-phase null-QS corresponds to a similar integral identity as before, but without any measure
$$
\int_{\partial \Omega^+} h(x) \ d\sigma_x - \int_{\partial \Omega^-} h(x) \ d\sigma_x = 0,
$$
where $\Omega^+ \cap \Omega^- =\emptyset.$

 As shown in \cite{ASW} the two-plane solutions (as they called it) are given by
\begin{enumerate}
\item\label{item1} $u=x_n^+$,
\item\label{item2} $u=-x_n^-$,
\item\label{item3} $u=x_n^+ -  (x_n+\gamma)^-$  for some real number $\gamma>0$,
\item\label{item4} $u=a x$.
\end{enumerate}
All these are global minimizers; for the last example to be a minimizer one needs   $a\ge 1$ (see \cite{ASW}, Lemma 4).
 It is, however, not clear whether these are the only two-phase null QS. Indeed, a (null)-QS does not need to be a minimizer of our functional, and the  only  requirement is that   it  satisfies a quadrature identity.


\section{Multi phase QS}
\subsection{The model equation}
It is apparent that once the  seed of the idea of two-phase QS (or any free boundary problem) has taken root, one 
may think of more complex situations where multi-phases are involved. In this section we shall rely on the above results for two-phase QS case, and provide  setting of a multi-phase problem, as done previously in segregation problems \cite{CTV}, 
or quadrature domain theory \cite{AvetHenrik2016}. Recently, the same approaches have been applied in shape optimization problems as well \cite{bucurmulti}.

The exact formulation of the multi-phase problem is done using the two phase version of it  as follows: \newline
Given  $m$ positive measures $\mu_i$, we want to find functions $u_i\geq 0, (i=1,\dots,m),$ with mutually disjoint supports $\Omega_i=\{u_i>0\}$ such that $\supp(\mu_i)\subset\Omega_i$ and
  \begin{align}\label{multi-QS-eq}
  \begin{cases}
  -\Delta(u_i-u_j)=\mu_i-\mu_j
 \;\;&\mbox{in}\;\;\R^N\setminus \overline{\cup_{k\neq i,j}\Omega_k},\\
|\nabla u_i  | = g, \;\; &\mbox{on}\;\; \partial  \Omega_i  \setminus  \overline{\cup_{k\neq i,j}\Omega_k}.
\end{cases}
\end{align}
In other words, for   each pair $(i,j)$ with $i\neq j$, the function
 $u_i - u_j$ solves a two-phase versions of our problem outside the union of the supports of the other functions.
 
A natural question that arises is: does the proposed model cover the two phase case? The answer is yes, because as for Multi-phase Quadrature domains (Theorem $2$ in \cite{AvetHenrik2016}), one can show the similar equivalence result.
 
\subsection{Existence of minimizers for multi-phase case}

In this section we will adapt the existence analysis, which has been done for multi-phase quadrature domains \cite{AvetHenrik2016}.
We start with the definition of  the minimization sets $K$ and $S.$ Define
\[
K=\{(u_1,u_2,\dots,u_m)\in (H^1(\R^N))^m \ :   \;\; u_i\geq 0,\; \mbox{for all}\;\;
i=1,\cdots ,m \},
\]
and
\[
S=\{(u_1,u_2,\dots,u_m)\in (H^1(\R^N))^m\ : \;\; u_i\geq 0,\;\mbox{and}\;\; u_i\cdot u_j=0,\; \mbox{for all}\;\; i\neq j  \}.
\]
Obviously  we have $S\subset K$. Next we define
\begin{equation}\label{mainfunctinal}
G(u_1,\dots,u_m)=\sum_{i=1}^m\int_{\R^N}\left(|\nabla u_i|^2 - 2f_i\cdot u_i+g^2 \chi_{\{u_i > 0\}}\right)dx,
\end{equation}
where each function $f_i$ and $g$ are satisfying conditions $(\tilde{A1})-(\tilde{A4}).$\footnote{For a general case of this functional one needs to replace $g^2 \chi_{\{u_i > 0\}}$ with   $\sum_{i=1}^mg_i^2 \chi_{\{u_i > 0\}}$ for $g_i$, for appropriate choices of $g_i$.}
 \be\left.\begin{array}{lll}
 	(\tilde{A1})&& f_i, g \in L^\infty(\R^N) ~{\rm ~for ~all ~i=1,2,...,m} \\
 	(\tilde{A2})&& {\rm supp\/} ~f_i^+~{\rm~ is ~compact ~for ~all ~i=1,2,...,m} \\
 	(\tilde{A3})&& g \geq 0\\
 	(\tilde{A4})&& {\rm either ~for ~all ~i~we ~have} ~f_i \leq - c_i< 0{~\rm  \ or  \ \/}~~~ g \geq  c_0 > 0~\\
 	&&~~~~~~~{\rm hold~outside~a~compact~set~for ~ some~positive~constants\/~}c_0,~c_i.
 \end{array} \right\}\ee
 	
 In light of Lemma $1$ in \cite{AvetHenrik2016}, one can show that for every  minimizer $(u_1,\dots,u_m)$ of  $G$ over $K,$ each component $u_i$ is going to minimize corresponding one-phase functional 
  \be
 J^1_{\{f_i,g\}}(\vp) =   \int\limits_{\R^N} \left( |\nabla \vp|^2 - 2f_i \vp  + g^2 \chi_{\{\vp > 0\}} \right)  \, dx,
 \ee 
 over the set $\{\vp \in H^1(\R^N) : \vp \geq 0 \}.$

 Hence, following \cite{gs} we say that the vector $(u_1,\dots,u_m)$ is a largest (smallest) minimizer of $G$ over $K,$ if for every $i,$ each component $u_i$ is accordingly the largest (smallest) minimizer (in the sense considered in \cite{gs}) of $J^1_{\{f_i,g\}}(\vp)$ over  the set $\{\vp \in H^1(\R^N) : \vp \geq 0 \}.$ 
 	
\begin{theorem}\label{embeddtheorem}
	Let $f_i(x),g(x)$ satisfy the conditions $(\tilde{A1})-(\tilde{A4}).$ Then $G(u_1,\dots,u_m)$ has at least one minimizer
	$(\bar{u}_1,\bar{u}_2,\dots,\bar{u}_m)$ in $S$, and also  all minimizers have compact support. 
	Moreover, the following  inclusion of supports  holds: 
	For any minimizer $(\bar{u}_1,\bar{u}_2,\dots,\bar{u}_m)$ of $G$  over $S$, and 
	the largest minimizer $(v_1,v_2,\dots,v_m)$ of $G$ over  $K$, we have 
	\begin{equation}\label{support}
	supp(\bar{u}_i)\subseteq supp(v_i), \:\:  i=1,\cdots,m.
	\end{equation}
	\end{theorem}
	\begin{proof}
The functional $G(u_1,u_2,\dots,u_m)$ is lower semi-continuous, coercive and convex. Since the set $S$ is closed, then the existence of a minimizer follows from  standard arguments of calculus of variations. Note that the minimizer is not necessarily unique.
		For simplicity, we make the following notations:
		\begin{align*}
		\bar{U}&\equiv(\bar{u}_1,\bar{u}_2,\dots,\bar{u}_m),\;\;V\equiv(v_1,v_2,\dots,v_m),\\
		\min(\bar{U},V)&\equiv(\min(\bar{u}_1,v_1),\min(\bar{u}_2,v_2),\dots,\min(\bar{u}_m,v_m)),\\
		\max(\bar{U},V)&\equiv(\max(\bar{u}_1,v_1),\max(\bar{u}_2,v_2),\dots,\max(\bar{u}_m,v_m)).
		\end{align*}
		To  see the  ordering  of the supports (equation \eqref{support}) one can proceed as in Lemma \ref{l1},
		 which clearly will imply the following inquality
		\[
		G(\min(\bar{U},V))+G(\max(\bar{U},V))\leq G(\bar{U})+G(V).
		\]
		Since $\bar{U}\in S$ and $V\in K,$ then $\min(\bar{U},V)\in S.$ Therefore $$G(\min(\bar{U},V))\geq G(\bar{U}),$$ which implies
		$$
		G(\max(\bar{U},V))\leq G(V).
		$$
		Observe that $\max(\bar{U},V)\in K$ and $V=(v_1,v_2,\dots,v_m)$ is the largest minimizer to $G(u_1,u_2,\dots,u_m)$ in $K.$ Hence,
		$$
		\max(\bar{U},V)\leq V,
		$$
		which is equivalent to
		\[
		\max(\bar{u}_i,v_i)\leq v_i,
		\]
		for all $i=1, \cdots , m.$ Thus $\bar{u}_i\leq v_i,$ which leads to
		\[
		supp(\bar{u}_i)\subseteq supp(v_i)
		\]
		for all $i=1, \cdots , m.$
		We recall that each component $v_i$ is the largest  minimizer to the functional $J^1_{\{f_i,g\}}(\vp)$ over  the set $\{\vp \in H^1(\R^N) : \vp \geq 0 \}.$  For these  functionals and  under more general setting  it has been  proved (see \cite[Theorem 1.4]{gs} ) that all minimizers have support in a fixed compact set. Thus, $supp(v_i)$ is compact, which in turn yields the compactness of $supp(\bar{u}_i),$ for all $i=1, \cdots , m.$  This completes the proof of  Theorem.
	\end{proof}
	
	Following the proof of Proposition $1$ in \cite{AvetHenrik2016}, we can prove a similar result for the functional \eqref{mainfunctinal}.
	\begin{proposition}\label{localQS}
		If $(\bar{u}_1,\bar{u}_2,\dots,\bar{u}_m)$  is a minimizer to the functional \eqref{mainfunctinal} subject to the set $S,$ then the following holds in the sense of distributions:
		\begin{align}\label{system1}
			\begin{cases}
		\Delta(\bar{u}_i-\bar{u}_j)=-f_i\chi_{\{\bar{u}_i>0\}}+ f_j\chi_{\{\bar{u}_j>0\}}\;\;\qquad & \mbox{in}\;\;\R^N\setminus  \cup_{k\neq i,j}\overline{\Omega}_k,\\
		|\nabla \bar{u}_i| = g, \qquad & \hbox{ on } \partial  \Omega_i  \setminus \cup_{k\neq i,j} \overline{\Omega_k},
		   \end{cases}
		\end{align}
		where $\Omega_i=\{\bar{u}_i>0\}.$ 
	\end{proposition}
	
	In the light of Two-phase QS and Proposition \ref{localQS}, we give a definition of Multi-phase version as follows 
	 \begin{definition}[Multi-phase Quadrature surface]
	 	Suppose we are given $m$ bounded positive measures $\mu_i$ and disjoint domains $\Omega_i$ such that $supp(\mu_i)\subset \Omega_i.$ If for  every harmonic functions $h\in HL^1(\overline{\Omega_i\cup\Omega_j}),$ such that $h$ is continuous across $\partial \Omega_i \cap \partial\Omega_j$,  
	 	and $h=0$ on $ \cup_{k\neq i, j} \partial \Omega_k,$ the following QI holds
	 	\begin{equation}\label{Q-identity-gen}	
	 	\int_{\partial \Omega_i} g h (x) \  d\sigma_x - \int_{\partial \Omega_j} g h (x) \  d\sigma_x= \int h d(\mu_i-\mu_j) \; ,
	 	\end{equation}
	 	then we call $\{\partial\Omega_i\}_{i=1}^m$ an  {\it m-phase QS} with respect to the measure $\{\mu_i\}_{i=1}^m,$ and a given  smooth positive function $g.$
	 	
	 	If we extend  the  test class $h$ to  the subharmonics in $\Omega_i$ and super-harmonics in $\Omega_j$ (due to negative sign in front of the integral) then  the equality  in \eqref{Q-identity-gen} is replaced with an inequality $(\geq) $.
	 \end{definition}
	
	The analogue of the Theorem \ref{existenceQS} for multi-phase case is the following result below.

	\begin{theorem}\label{existQD-gen}
		Let $\mu_i$ be  given Radon measures with compact supports,  that satisfy Sakai's condition as in Definition \ref{sakai}. Suppose that for each $\mu_i$ the corresponding  one-phase quadrature surface $\partial Q_i$ (see \eqref{1-Q-identity})  is such that
		\begin{equation}
		\overline{Q}_i\cap supp(\mu_j)=\emptyset,\;\; \mbox{for every}\;\; i\neq j.
		\end{equation}
		Then, we have a solution to our multi-phase free boundary problem \eqref{multi-QS-eq} along with  $supp(\mu_i)\subset supp(\{u_i>0\}),$ for all $l=1,2\dots,m.$
	\end{theorem}
	
	The proof of this result repeats the same lines as in Theorem \ref{existenceQS}, and therefore is omitted.

\subsection{Analysis of junction points }
In this section our goal is to show the absence of triple junction points in $\R^N$, away from  the support of the measures $\mu_i.$ More exactly we shall show that for multi-phase QS, there is at most two phases that can meet at each point. 
 In the case of multi-phase quadrature domains it was shown (see \cite{AvetHenrik2016})  that a triple junction may actually appear.
 
For the proof of the main result,  in this section, we will need the multi-phase counterpart (see \cite{multiphase-note}) of  a celebrated Caffarelli-Jerison-Kenig monotonicity formula \cite{CJK}. 
 
\begin{theorem}(\cite{multiphase-note})[Three-phase monotonicity formula]\label{monotonicity-ACF}
Let $B_1\subset \R^N$ be the unit ball in $ \R^N$ and let $ u_i\in H^1(B_1), i = 1,2, 3, $ be three non-negative Sobolev functions such that
\[ \Delta u_i+1\geq 0, \qquad \forall i=1,2,3, \qquad \mbox{and} \qquad u_i\cdot u_j = 0\;\; \mbox{a.e.}\;\;  in\; B_1,\;\;  \forall i\neq j.\]
Then there are dimensional constants $\varepsilon> 0$  and $C_N > 0$  such that for each $r \in( 0, 1)$  we have
\[ \prod_{i=1}^3 \left( \frac{1}{r^{2+\varepsilon}}\int_{B_r}\frac{|\nabla u_i|^2}{|x|^{N-2}}dx \right) \leq C_N\left(1+ \sum_{i=1}^3\int_{B_1}\frac{|\nabla u_i|^2}{|x|^{N-2}}dx \right)^3. 
\]

\end{theorem}

\begin{lemma}[Non-degeneracy]\label{non-deg}
	Let  $(u_1,u_2,\dots\,u_m)\in S $ be a minimizer to \eqref{mainfunctinal}.
	Then there exist a  constant $D_{N,f_i,g}>0,$  depending on dimension $N,$ and functions $f_i,g,$ such that for every $x_i\in\partial\Omega_i\cap B_{\frac{1}{2}}(0)$ we have
	
	$$
	\oint\limits_{\partial B_{r}(x_i)} u_i \geq r\cdot D_{N,f_i,g}.
	$$
	Here,  we set  $\Omega_i=\{u_i>0\},$ and $i=1,2,\dots, m.$ 

\end{lemma}
\begin{proof}

	To see this for some fixed $i$,  we set  
	$$
	G_{r,i}(v)=\int_{B_{r/2}}\left(|\nabla v|^2-2f_iv+ g^2 \chi_{\{v>0\}}\right) dx,
	$$
	and
	$$
	\tilde{G}_{r,i}(v)=\int_{B_{r/2}}\left(|\nabla v|^2-2M_iv+ l^2 \chi_{\{v>0\}}\right) dx,
	$$
	where $l=\underset{{B_r}}{\inf}\;g,\; M_i=\underset{{B_r}}{\sup}\; {f_i^+}.$ 

	Now for a constant $\beta>0$, if we define  
	$$
	K_{\beta,i}=\{v\in H^1(B_{r/2}): \quad  v\geq 0, \quad  v=\beta\;\;\mbox{on}\;\;\partial B_{r/2} \},
	$$
	then following the proof of Lemma $2.8$ in \cite{gs}, we conclude that the largest minimizer   
	$ v_{\beta,i}$,	of $	\tilde{G}_{r,i}$ over $K_{\beta,i},$   vanishes  in the ball $B_{r/4},$ provided $r$ and $\beta$ are small enough. The upper thresholds for the constants $r$ and $\beta$  can be taken as follows:
	$$
	0<r<2\frac{Nl}{\hat{M}}\;\;\mbox{and}\;\;0<\beta\leq \beta_0(r,l,M_i),
	$$ 
	where $\hat{M}=\underset{i}{max}\;M_i $ and  $\beta_0(r,l,M_i)$ is a corresponding threshold of one  phase problem with force term $f_i$ (see the proof of Lemma $2.8$  in \cite{gs}). On the other hand due to Harnack inequality for the component $u_i$ we have 
	$$
	u_i\leq C_1\oint\limits_{\partial B_{r/2}} u_i+C_2r^2M_i,
	$$
	where the constants $C_1$ and $C_2$ are depending only on dimension $N.$  Using the following rescaling property $\beta_0(r,l,M_i)=r\beta_0(1,l,rM_i),$ one can easily achieve 
	\begin{equation}\label{threshold}
	u_i(x)<\beta_0(r,l,M_i)\;\;\mbox{on}\;\;\partial B_{r/2},
	\end{equation}
	by letting 
	$$
	C_1 \frac{1}{r}\oint\limits_{\partial B_{r}} u_i+rC_2M_i<\beta_0(1,l,rM_i).
	$$
	Since $\beta_0(1,l,0)>0,$ then taking $r$ small enough we will have \eqref{threshold}, provided $ \frac{1}{r}\oint\limits_{\partial B_{r}} u_i\leq C,$ where $C>0$ is some  constant. Our aim is to prove that $u_i=0$ in $B_{r/4}.$  We take $v_{\beta,i}$ to be the largest minimizer of $\tilde{G}_{r,i}$ over the set $K_{\beta,i}$  for $\beta=\beta_0(r,l,M_i).$ We define a new function $w$ to be $\min(u_i,v_{\beta,i})$ in $B_{r/2}$ and equal $u_i$ in $\R^N\setminus B_{r/2}.$ The inequality \eqref{threshold} implies $(u_1,\dots,u_{i-1},w,u_{i+1},\dots,u_m)\in S,$ and therefore
	$$
	G(u_1,\dots,u_m)\leq G(u_1,\dots,u_{i-1},w,u_{i+1},\dots,u_m).
	$$
	This leads to
	\begin{equation}\label{ineq1}
	G_{r,i}(u_i)\leq G_{r,i}(\min(u_i,v_{\beta,i})).
	\end{equation}
	
	Using the same arguments as in the proof of Lemma \ref{l1}, we can obtain the following inequality
	\begin{equation}\label{ineq2}
	G_{r,i}(\min(u_i,v_{\beta,i}))+\tilde{G}_{r,i}(\max(u_i,v_{\beta,i}))\leq G_{r,i}(u_i)+\tilde{G}_{r,i}(v_{\beta,i})
	\end{equation}
	
	Thus, in the light of \eqref{ineq1} and \eqref{ineq2} we get $\tilde{G}_{r,i}(v_{\beta,i})\geq \tilde{G}_{r,i}(\max(u_i,v_{\beta,i})),$ which in turn implies $\max(u_i,v_{\beta,i})\leq v_{\beta,i}.$ The latter inequality follows from the fact that $v_{\beta,i}$ is a largest minimizer to $\tilde{G}_{r,i}$ over the set $K_{\beta,i},$ and $\max(u_i,v_{\beta,i})\in K_{\beta,i}.$ Hence, $u_i\leq v_{\beta,i}$ in $B_{r/2}$ and this gives that $u_i=0$ in $B_{r/4}.$ Thus, we have proved that for every component $u_i$  there exists a dimensional  constant $C_N>0,$ depending also on $l$ and $M_i$  such that for every sufficiently small $r>0$  the following statement is true:
	\begin{equation}\label{non-deg-cond}
	\frac{1}{r}\oint\limits_{\partial B_{r}} u_i\leq C_N\; \Rightarrow\; u_i=0\;\; \mbox{in}\;\; B_{r/4}.
	\end{equation}
	This basically gives the desired non-degeneracy property. Similarly, it can be shown that the statement \eqref{non-deg-cond} remains true  with $B_{kr}$ in place of  $B_{r/4}$ for any $0<k<1.$ In this case the constant $C_N$ will also depend on $k.$  
\end{proof}

\begin{lemma}[Lipschitz regularity]\label{Lip-regular}
	Let $U=(u_1,u_2,\dots\,u_m)\in S $ be a minimizer to \eqref{mainfunctinal}, and suppose 
	$\hbox{supp} (\mu) \cap B_2 (0) = \emptyset$. 
	Then there exists a universal constant $C_N >0$ depending only on $N,$ such that
	$$
	||u_i||_{C^{0,1}(B_{1/4})}\leq C_N||U||_{L^{2}(B_{1})},
	$$
	for every $i=1,2,\dots,m.$
\end{lemma}
\begin{proof}

	  Invoking Theorem $7.1$ 	  in \cite{fbp-localization}, by choosing  $F(W)$ as in equation $(1.6)$ in the same paper, we may conclude that $u_i$ is locally $C^\alpha$. From here one may apply Lemma $5.2$ in \cite{ACF} to conclude 	  
	 	\begin{equation}\label{mean-value}
	 	\int_{\partial B_r(z)} |u_i| \leq Cr^N, \quad z\in \Gamma_i\cap\Gamma_j,
	 	\end{equation}
	 	with $C$ universal constant, depending only on the distance between $z$ and the support of measures (in our case).  It should be remarked that in Lemma 5.2 of \cite{ACF} we have to take only two functions at a time, so as to apply the monotonicity function. The latter can be found in more advanced forms in \cite{CJK}. Further the importance of initial Holder regularity is needed in Lemma 5.2 (equation (5.12)) in \cite{ACF}.

	 	From \eqref{mean-value} we may now infer Lipschitz regularity as done in the proof of Theorem $5.3$ in \cite{ACF}, where one also needs that $|u_i|$ is a sub solution, which is fulfilled by our solutions.

\end{proof}

\begin{remark}
It is noteworthy that several recent papers prove Lipschitz regularity of solutions for two and multi-phase problems with heuristic arguments, without stressing the need for initial partial regularity.
 It needs to be stressed that the conditions on the functions in  the monotonicity formula of \cite{ACF}, and the succeeding ones, have been relaxed  considerably, and in general one can avoid continuity of solutions. Nevertheless, for applying the formula to prove regularity of solutions in free boundary problems, one does need to begin with some initial partial regularity.  This part of the problem is too often neglected and not taken seriously.  
This has been highlighted in  our proof of Lipschitz regularity of solutions in \cite{fbp-localization} where we begin with solutions satisfying H\"older regularity. To the author's best knowledge, 
 the $C^\alpha$-regularity for multi-phase problems is by no means an easy problem, and cannot be done  as that of  Theorem $2.1$ in \cite{ACF}.
\end{remark}

\begin{theorem}
 Let $(u_1,u_2,\dots\,u_m)\in S $ be a minimizer to \eqref{mainfunctinal}.
 Then there is a universal constant $R_0 >0 $ (depending only on the norms) such that for any  point $z_{i,j} \in \partial \{u_i >0\} \cap \partial \{u_j >0\}$
we have $|z_{i_1i_2} -z_{i_3i_4}| > R_0$, provided $(i_1,i_2) \neq (i_3,i_4)$. Here $i_k \in \{1, \cdots, m\}$. 

In particular  triple junction points cannot appear, and that  two different class of two-phase points stay uniformly away from each other.
\end{theorem}

This theorem can be set in relation  to  segregation problems that have been in focus lately, see \cite{bucurmulti}.  A particular application of this theorem is that in segregation problems, where  multi-phase Bernoulli type free boundaries appear in the limit, one can claim that more than two phases cannot meet at the same time.

\begin{proof}

We first notice that by compactness, and non-degeneracy, if there is a  sequence $z^k_{i_1,i_2}, $ 
$z^k_{i_3,i_4},  $  ($k=1, 2, \cdots $) of  two-phase points of different classes such that 
$|z^k_{i_1,i_2} - z^k_{i_3,i_4} | \to 0$, then the limit  point $w= \lim_k z^k_{i_1,i_2} =  \lim_k z^k_{i_3,i_4}$ is a triple junction point. Hence to prove the theorem, it suffices to show that  triple junction points do not exist. 
	 
By non-degeneracy (see Lemma \ref{non-deg}) for  each $z_i\in\partial\{u_i>0\}$ 
there exists a point $y_i\in\partial B_{r/2}$ with $u_i(y_i)\geq D_{N,f_i,g}\cdot r.$ 

From Lemma 7.2, $u_i$ is Lipschitz regular, and therefore  
 $u_i(x)>0$ in $B_{\varepsilon r}(y_i),$ for some small enough $\varepsilon>0.$ Thus  there exists a positive constant $c_0>0$ such that for any small enough $r>0$ the following property holds:
 $$
 \cfrac{|\{u_i>0\}\cap B_r(z_i)|}{|B_r(z_i)|} \geq c_0>0.
 $$
 Since all sets $\{u_i>0\}$ are disjoint, then there exists  a dimensional constant $\alpha_0>0$ 
 $$
 \cfrac{|\{u_i=0\}\cap B_r(z_i)|}{|B_r(z_i)|} \geq \alpha_0>0,
 $$
 for  every $i=1,2,\dots,m.$

 We will need the following version of Poincare inequality: For every function $v\in H^1(B_r)$ we have
 \begin{equation}\label{Poincare-gen}
 |\{v=0\}\cap B_r|\left(\frac{1}{r}\oint\limits_{\partial B_{r}} v\right)^2\leq C_N \int\limits_{B_r}|\nabla v|^2.
 \end{equation}
 The proof of this inequality can be found implicitly in \cite[Lemma $3.2$]{AC}. Another reference is Lemma $4.5$ in \cite{bucurmulti}. 
 In view of non-degeneracy property and inequality \eqref{Poincare-gen} we arrive at:
 $$
 \alpha_0\cdot D_{N,f_i,g}^2\leq \cfrac{|\{u_i=0\}\cap B_r(z_i)|}{|B_r(z_i)|}\left(\frac{1}{r}\oint\limits_{\partial B_{r}} u_i\right)^2\leq \cfrac{C_N}{|B_r(z_i)|}\int\limits_{B_r}|\nabla u_i|^2.
 $$
 Thus, there exists a universal constant $L_i>0$ depending only on  $N,f_i,g,$ such that  
 
 \begin{equation}\label{gradient-below} 	
\int\limits_{B_r(z_i)}|\nabla u_i|^2\geq L_i\cdot r^N.
\end{equation}
 
Now,  let the origin be a possible triple junction point for components $u_{i_1},u_{i_2}$ and $u_{i_3},$ away from the measures $\mu_i.$ Our aim is to apply the multi-phase version of Caffarelli-Jerison-Kenig monotonicity formula around the origin and come to a contradiction.  Since the triple junction point is away  from the measures $\mu_i$, the constants $L_i$ do not depend  on $f_i$ in a small neighborhood of the origin and therefore are the same. 

First we recall  the following inequality obtained in \cite[Remark $1.5$]{CJK} (see also  \cite{multiphase-note}):
 Suppose that $u\in H^1(B_2)$ is a nonnegative Sobolev function such that $\Delta u + 1 \ge 0$ on $B_2 \subset\R^N.$ Then, there is a dimensional constant $Q_N>0,$ such that
\begin{equation}\label{aux_lemma-CJK}
\int\limits_{B_1}{\frac{|\nabla u|^2}{|x|^{N-2}}}dx\leq Q_N\left(1+ \int\limits_{B_2} u^2dx\right).
\end{equation}
 According to Theorem \ref{monotonicity-ACF},  estimate \eqref{gradient-below} and inequality \eqref{aux_lemma-CJK} we obtain
\begin{align*}
r^{-3\varepsilon}L_1^3&\leq   \prod_{j=1}^3 \left( \frac{1}{r^{N+\varepsilon}}\int\limits_{B_r(0)}{|\nabla u_{i_j}|^2}dx \right)\leq \prod_{j=1}^3 \left( \frac{1}{r^{2+\varepsilon}}\int\limits_{B_r(0)}{\frac{|\nabla u_{i_j}|^2}{|x|^{N-2}}}dx \right)\leq\\& \leq C_N\left(1+ \sum_{j=1}^3\int\limits_{B_1(0)}{\frac{|\nabla u_{i_j}|^2}{|x|^{N-2}}}dx \right)^3\leq C_N\left((1+3Q_N)+Q_N\sum_{j=1}^{3}\int\limits_{B_2(0)} u_{i_j}^2dx\right)^3. 
\end{align*}

By letting   $r\to 0^+,$ we  conclude  that $\sum_{j=1}^{3}\int\limits_{B_2(0)} u_{i_j}^2dx=+\infty,$ which gives a contradiction and completes the proof.

\end{proof}
\begin{remark}
Note that without loss of generality, we have assumed $ \hbox{supp} (\mu_i) \cap B_2(0) = 
\emptyset$ and therefore we used $\Delta u_i\ge -1$ condition to obtain the main result.
This assumption maybe justified by scaling, once the free boundary is a certain distance $r_0 >0$
away from the support of the measures.
\end{remark}


\bibliographystyle{acm}

\bibliography{multi-QS}

\end{document}